*Research Article*

# Variations on Hammersley's interacting particle process

Arda Atalik[1], H. S. Melihcan Erol[2], Gökhan Yıldırım[3,*], Mustafa Yilmaz[1]

[1]*Department of Electrical and Electronics Engineering, Bilkent University, Ankara, TR 06800, Turkey*
[2]*Electrical Engineering & Computer Science, MIT, Cambridge, MA 02139, USA*
[3]*Department of Mathematics, Bilkent University, Ankara, TR 06800, Turkey*

**Abstract**

The longest increasing subsequence problem for permutations has been studied extensively in the last fifty years. The interpretation of the longest increasing subsequence as the longest 21-avoiding subsequence in the context of permutation patterns leads to many interesting research directions. We introduce and study the statistical properties of Hammersley-type interacting particle processes related to these generalizations and explore the finer structures of their distributions. We also propose three different interacting particle systems in the plane analogous to the Hammersley process in one dimension and obtain estimates for the asymptotic orders of the mean and variance of the number of particles in the systems.

## 1. Introduction

The longest increasing subsequence problem for permutations, also called Ulam's problem [18], has motivated an interesting research program at the intersection of different branches of mathematics such as probability theory, random matrix theory, operator theory, and statistical physics [2, 5, 9, 14, 16]. Any arrangements of the elements in $[n] := \{1, 2, \cdots, n\}$ is called a permutation that can also be considered a one-to-one and onto function from $[n]$ to itself. We use $S_n$ to denote the set of all permutations of length $n$. For $\sigma = \sigma_1\sigma_2\cdots\sigma_n \in S_n$, we define $\ell is_n(\sigma)$ as the length of the longest increasing subsequence in $\sigma$, that is, the maximum $k \in [n]$ for which the conditions $1 \leq i_1 < i_2 < \cdots < i_k \leq n$ and $\sigma_{i_1} < \sigma_{i_2} < \cdots < \sigma_{i_k}$ hold. A well-known result in combinatorics, called Erdős-Szekeres's lemma, states that any sequence of $n^2 + 1$ distinct numbers has either an increasing or decreasing subsequence of length at least $n + 1$. It follows from this fact that the expected value of $\ell is_n$, $\mathbf{E}(\ell is_n)$, grows asymptotically proportional to $c\sqrt{n}$ for some positive constant $c$ as $n$ tends to infinity. In a seminal paper in 1972 [7], Hammersley introduced a Poissonized version of the problem and proved the existence of the limit for the sequence $\mathbf{E}(\ell is_n)/\sqrt{n}$; later it was shown that $c = 2$, for the details, see [2]. A complete solution for the problem was provided in 1999 [3] by determining the order of the fluctuations around the mean as $n^{1/6}$ and limiting distribution as the Tracy-Widom GUE distribution, that is,

$$\lim_{n\to\infty} \mathbf{P}\left(\frac{\ell is_n - 2\sqrt{n}}{n^{1/6}} \leq x\right) = F_2(x)$$

for all $x$. The Tracy-Widom distribution first appeared as the limiting distribution for the rescaled largest eigenvalue of the Gaussian unitary ensemble (GUE) from the random matrix theory [17]. Other Tracy-Widom distributions from orthogonal (GOE) and symplectic (GSE) matrix ensembles also appear as a limiting distribution for some specific permutations classes, see the introduction of [4].

It also follows from the Hammersley's work [2, 7] that a simple interacting particle process on the unit interval in which the macroscopic quantity defined as the number of particles in the system has the same statistical distribution with the random variable $\ell is_n$. The particle process approach also gives a very efficient and elegant algorithm for simulating $\ell is_n$ and described as follows: Initially, there are zero particles in the system. At each step, a particle appears at a uniform random point $u$ in the interval $[0, 1]$; simultaneously the nearest particle (if any) to the right of $u$ disappears. If $p_n^h$ denotes the number of particles in the system after $n$ steps, then $p_n^h$ and $\ell is_n$ have the same probability distribution hence the large time behavior of the particle system follows the Tracy-Widom GUE distribution.

---

*Corresponding author (gokhan.yildirim@bilkent.edu.tr).

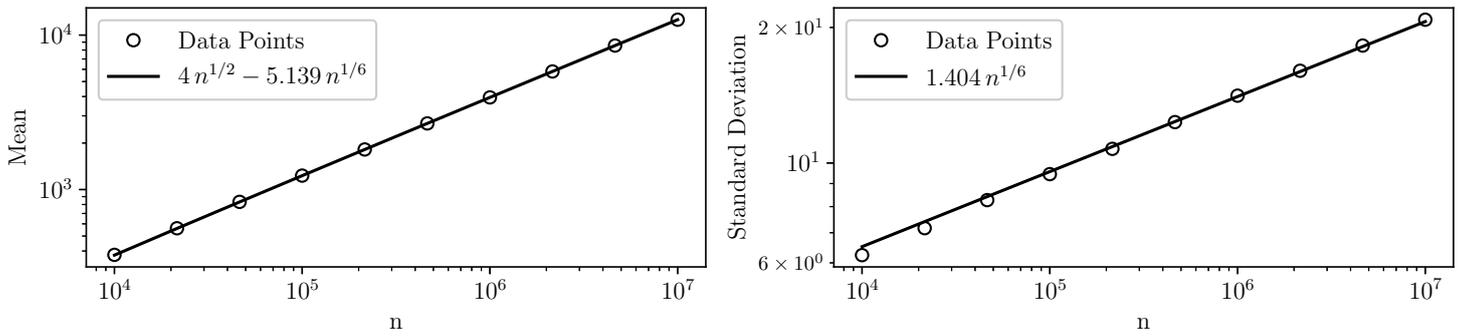

Figure 1: Log-Log plots of the best fits for the mean and standard deviation of the number of particles in two-aligned segments. The data indicate that $E(p_n^I) \approx 4 n^{1/2} + a n^{1/6}$ and $SD(p_n^I) \approx b n^{1/6}$ with 95% confidence intervals for $a \in [-5.222, -5.057]$ and $b \in [1.392, 1.417]$ with corresponding MSEs 1.054 and 0.024, respectively.

|  | $n_1$ | $n_2$ | $n_3$ | $n_4$ | $n_5$ | $n_6$ | $n_7$ | $n_8$ | $n_9$ | $n_{10}$ |
|---|---|---|---|---|---|---|---|---|---|---|
| $(4n^{\frac{1}{2}} - \text{mean}(p_n^{321})) n^{-\frac{1}{6}}$ | 4.790 | 4.869 | 4.940 | 4.997 | 5.056 | 5.098 | 5.147 | 5.175 | 5.201 | 5.240 |
| (st. dev. $(p_n^{321})) n^{-\frac{1}{6}}$ | 1.345 | 1.358 | 1.380 | 1.387 | 1.388 | 1.401 | 1.411 | 1.410 | 1.412 | 1.418 |
| skewness | 0.131 | 0.116 | 0.130 | 0.129 | 0.098 | 0.102 | 0.116 | 0.120 | 0.115 | 0.106 |
| kurtosis | 3.024 | 3.027 | 3.040 | 3.030 | 2.992 | 3.047 | 3.033 | 3.006 | 3.011 | 3.017 |

Table 1: Monte Carlo simulation results for two-aligned segments for 10 logarithmically spaced $n$ values from $10^4$ to $10^7$.

### 1.1. Pattern avoidance and generalization of the Ulam's problem

For $\tau = \tau_1 \tau_2 \cdots \tau_k \in S_k$ and $\sigma = \sigma_1 \sigma_2 \cdots \sigma_n \in S_n$, it is said that $\tau$ appears as a *pattern* in $\sigma$ if there exists a subset of indices $1 \leq i_1 < i_2 < \cdots < i_k \leq n$ such that $\sigma_{i_s} < \sigma_{i_t}$ if and only if $\tau_s < \tau_t$ for all $1 \leq s, t \leq k$. If $\tau$ does not appear as a pattern in $\sigma$, then $\sigma$ is called a $\tau$-*avoiding* permutation. For example, $132 \in S_3$ appears as a pattern in $246513$ because it has the subsequences $24 - - - 3$, $2 - 6 - - 3$, $2 - 65 - -$, $2 - -5 - 3$ or $-465 - -$. On the other hand, $4213 \in S_4$ does not appear as a pattern in $246513$. We denote by $S_n(\tau)$ the set of all $\tau$-avoiding permutations of length $n$. More generally, for a set $T$ of patterns, we use the notation $S_n(T) = \bigcap_{\tau \in T} S_n(\tau)$. It is known that for any pattern $\tau \in S_k$, $f(\tau) := \lim_{n \to \infty} |S_n(\tau)|^{1/n}$ exists as a finite number [12]. For some specific cases, it is possible to calculate the limiting value explicitly such as $f(\tau) = 4$ for any $\tau \in S_3$ and $f(12 \ldots k) = (k-1)^2$. Albert introduced an interesting generalization and reinterpretation of the longest increasing subsequence problem in the context of pattern avoidance [1]. For a given pattern $\tau \in S_k$ and $\sigma \in S_n$, let $\ell_n^\tau(\sigma)$ be the length of the longest $\tau$-avoiding subsequence in $\sigma$. Note that the longest increasing subsequence corresponds to the longest 21-avoiding subsequence, that is, $\ell_n^{21}(\sigma) = \ell is_n(\sigma)$. A set $T$ of patterns is called proper if $T$ does not contain both patterns $12 \ldots k$ and $j \ldots 21$ for some $k$ and $j$. The following conjecture was proposed in [1]:

**Conjecture:** Let $T$ be a proper set of patterns. Then

$$\lim_{n \to \infty} \frac{\mathbf{E}(\ell_n^T)}{\sqrt{n}} = 2\sqrt{\limsup_{n \to \infty} |S_n(T)|^{1/n}}.$$

This conjecture was proven in [1] for monotone patterns of length $k$, that is,

$$\lim_{n \to \infty} \frac{\mathbf{E}(\ell_n^{k(k-1) \cdots 21})}{\sqrt{n}} = 2(k-1).$$

In this note, we propose Hammersley-type particle processes corresponding to the pattern $\tau = k \ldots 21$, and numerically study the finer structures of the distribution of the rescaled random variables $\ell_n^\tau$ for $\tau = 4321$ and $\tau = 321$. Our particle processes give also an efficient algorithm for simulating $\ell_n^\tau$ for any monotone pattern $\tau$.

## 2. Numerical Results

### 2.1. Hammersley-type processes corresponding to the monotone pattern $\tau = k(k-1) \ldots 21$

Consider $k - 1$ unit intervals $[0, 1]$ in parallel and label them from top to bottom as $l_1, l_2, \ldots, l_{k-1}$. Initially, there are zero particles in the system. At each step, a particle appears at a uniform random point $u_o^1$ in the interval $l_1$; simultaneously the nearest particle $u_r^1$ (if any) to the right of $u_o^1$ disappears from $l_1$ and a particle appears in $l_2$ at $u_o^2 = u_r^1$; simultaneously the nearest particle $u_r^2$ (if any) to the right of $u_o^2$ disappears from $l_2$ and a particle appears in $l_3$ at $u_o^3 = u_r^2$; .... ; simultaneously the nearest particle $u_r^{k-1}$ (if any) to the right of $u_o^{k-1} = u_r^{k-2}$ disappears from $l_{k-1}$ and leaves the system. Let $p_n^\tau$ denotes the number of particles in the system after $n$ steps.

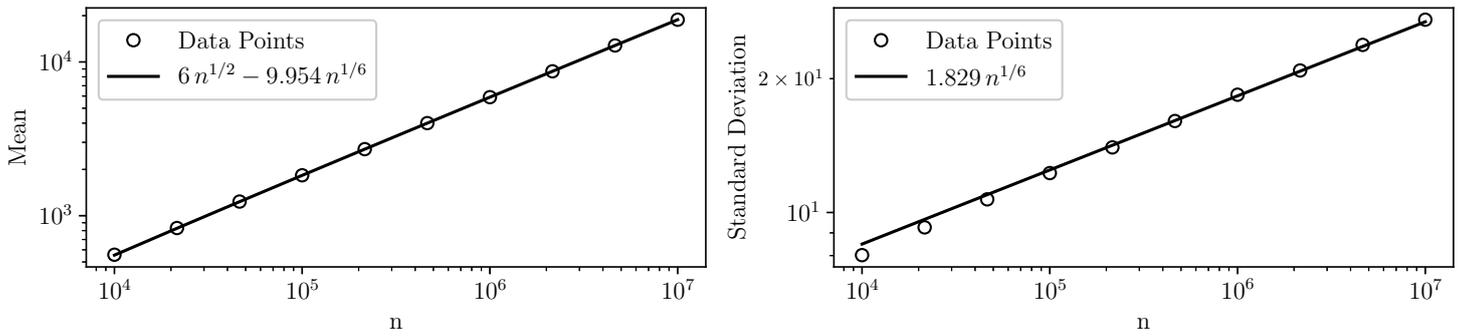

Figure 2: Log-Log plots of the best fits for the mean and standard deviation of the number of particles in three-aligned segments. The data indicate that $E(p_n^I) \approx 6\,n^{1/2} + a\,n^{1/6}$ and $SD(p_n^I) \approx b\,n^{1/6}$ with 95% confidence intervals for $a \in [-10.140, -9.774]$ and $b \in [1.809, 1.850]$ with corresponding MSEs 5.056 and 0.065, respectively.

|  | $n_1$ | $n_2$ | $n_3$ | $n_4$ | $n_5$ | $n_6$ | $n_7$ | $n_8$ | $n_9$ | $n_{10}$ |
|---|---|---|---|---|---|---|---|---|---|---|
| $(6n^{\frac{1}{2}} - \text{mean}(p_n^{4321}))\,n^{-\frac{1}{6}}$ | 9.188 | 9.363 | 9.519 | 9.647 | 9.767 | 9.866 | 9.963 | 10.034 | 10.093 | 10.172 |
| (st. dev. $(p_n^{4321}))\,n^{-\frac{1}{6}}$ | 1.727 | 1.755 | 1.786 | 1.800 | 1.810 | 1.825 | 1.840 | 1.836 | 1.844 | 1.849 |
| skewness | 0.087 | 0.077 | 0.084 | 0.089 | 0.069 | 0.072 | 0.078 | 0.084 | 0.079 | 0.074 |
| kurtosis | 2.998 | 3.002 | 3.020 | 3.020 | 2.979 | 3.039 | 3.027 | 3.005 | 3.006 | 3.020 |

Table 2: Monte Carlo simulation results for three-aligned segments for 10 logarithmically spaced $n$ values from $10^4$ to $10^7$.

**Theorem 2.1.** *Under the uniform probability distribution, the random variables $\ell_n^\tau$ and $p_n^\tau$ have the same distribution.*

*Proof.* The proof follows from the main result of [6], which is a generalization of Robinson-Schensted [15] correspondence, and the patience sorting algorithm [2]. Recall that a partition of integer $n$ is a collection of integers $\lambda = (\lambda_1, \lambda_2, \ldots, \lambda_l)$ such that $\lambda_1 \geq \lambda_2 \geq \ldots \geq \lambda_l \geq 1$ and $\sum_{i=1}^{l} \lambda_i = n$. A given partition $\lambda$ of $n$ can be represented by a Young diagram, a finite collection of boxes, or cells, arranged in left-justified rows, with the $i^{th}$ row having $\lambda_i$ boxes. A standard Young tableau is a filling of the diagram $\lambda$ by integers $1, 2, \cdots, n$ such that the numbers are in increasing order along each row and each column. The Robinson-Schensted correspondence uniquely associates a pair $(Q, R)$ of Young tableaux of the same shape of size $n$ to each $\sigma \in S_n$. Moreover, $lis_n(\sigma) = \lambda_1$. Then the result follows from Greene's theorem, $\ell_n^{k(k-1)\cdots 21}(\sigma) = \lambda_1 + \lambda_2 + \cdots + \lambda_{k-1}$, and the patience sorting algorithm. □

For the patterns $\tau = 4321$ and $\tau = 321$, we carried out Monte Carlo simulations for the corresponding particle systems for each of the 10 logarithmically spaced $n$ values from $n_1 = 10^4$ to $n_{10} = 10^7$. We obtained $10^5$ samples for each case and analyzed various statistical properties of $p_n^\tau$, the number of particles in the system after $n$ steps. We summarized the numerical results for $\tau = 321$, in Figure 1 and Table 1; for $\tau = 4321$, in Figure 2 and Table 2. The results are consistent with the verified cases of the conjecture for the order of the mean; they also indicate that the variance is of order $n^{1/6}$.

### 2.2. Hammersley-type processes in the plane

In this section, we introduce three Hammersley-type interacting particle processes in the plane and numerically study their statistical properties. We summarize our results in Tables 3, 4, and 5. We use $\mathcal{K}$ to denote the unit square in the Cartesian plane with corners at the points $(0, 0), (1, 0), (0, 1)$, and $(1, 1)$. The illustrations of the models are given in Figure 6.

**Model-I** For $(u, v) \in \mathcal{K}$ and a real number $m$, let $\mathcal{K}_m(u, v)$ denote the part of the unit square $\mathcal{K}$ above the line through $(u, v)$ with slope $m$. That is, $\mathcal{K}_m(u, v) = \{(x, y) \in \mathcal{K} : y > m(x - u) + v\}$. Fix a slope $m$. Initially, there are zero particles in the system. At each step, a particle appears at a uniform random point $(u, v)$ in the unit square $\mathcal{K}$; simultaneously the closest particle in $\mathcal{K}_m(u, v)$ to the point $(u, v)$ (if any) disappears. Let $p_n^I$ denote the number of particles in the system after $n$ steps.

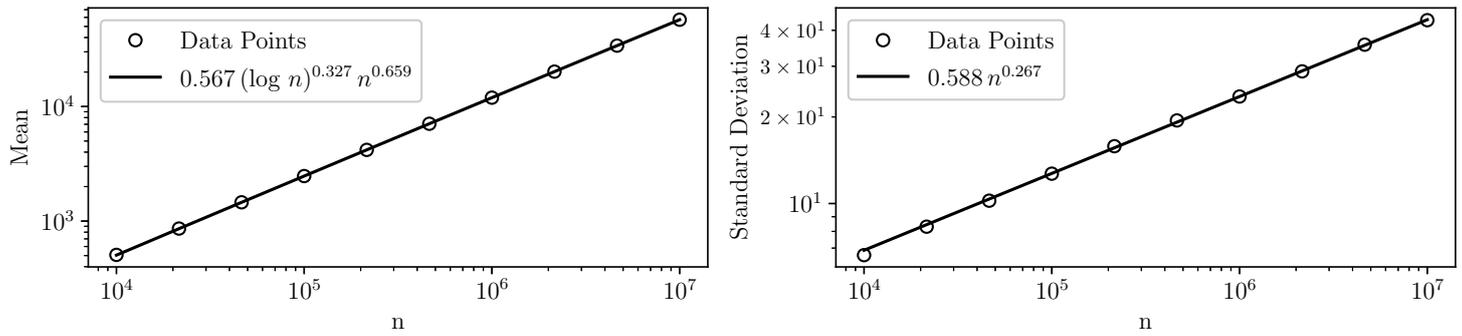

Figure 3: Log-Log plots of the best fits for the mean and standard deviation of the number of particles in model-I. The data indicate that $E(p_n^I) \approx a (\log n)^b n^c$ and $SD(p_n^I) \approx d n^f$ with 95% confidence intervals for $a \in [0.538, 0.595]$, $b \in [0.297, 0.357]$, $c \in [0.657, 0.661]$, and $d \in [0.555, 0.621]$, $f \in [0.263, 0.271]$ with corresponding MSEs 3.226 and 0.029, respectively.

|  | $n_1$ | $n_2$ | $n_3$ | $n_4$ | $n_5$ | $n_6$ | $n_7$ | $n_8$ | $n_9$ | $n_{10}$ |
|---|---|---|---|---|---|---|---|---|---|---|
| mean$(p_n^I) (\log n)^{-0.327} n^{-0.659}$ | 0.570 | 0.568 | 0.567 | 0.567 | 0.567 | 0.567 | 0.567 | 0.567 | 0.567 | 0.567 |
| (st. dev. $(p_n^I)) n^{-0.267}$ | 0.565 | 0.578 | 0.580 | 0.586 | 0.595 | 0.596 | 0.588 | 0.586 | 0.591 | 0.585 |
| skewness | 0.082 | 0.066 | 0.059 | 0.033 | 0.027 | 0.047 | 0.068 | 0.025 | 0.015 | 0.007 |
| kurtosis | 3.045 | 3.016 | 3.009 | 3.018 | 3.008 | 3.045 | 2.900 | 2.902 | 3.002 | 3.036 |

Table 3: Monte Carlo simulation results for the model-I.

**Model-II** For $(u, v) \in \mathcal{K}$, let $\mathcal{C}(u, v)$ denote the part of the unit square $\mathcal{K}$ outside the circle through $(u, v)$ with center $(0, 0)$. That is, $\mathcal{C}(u, v) = \{(x, y) \in \mathcal{K} : x^2 + y^2 > u^2 + v^2\}$. Initially, there are zero particles in the system. At each step, a particle appears at a uniform random point $(u, v)$ in the unit square $\mathcal{K}$; simultaneously the closest particle in $\mathcal{C}(u, v)$ to the point $(u, v)$ (if any) disappears. Let $p_n^{II}$ denote the number of particles in the system after $n$ steps.

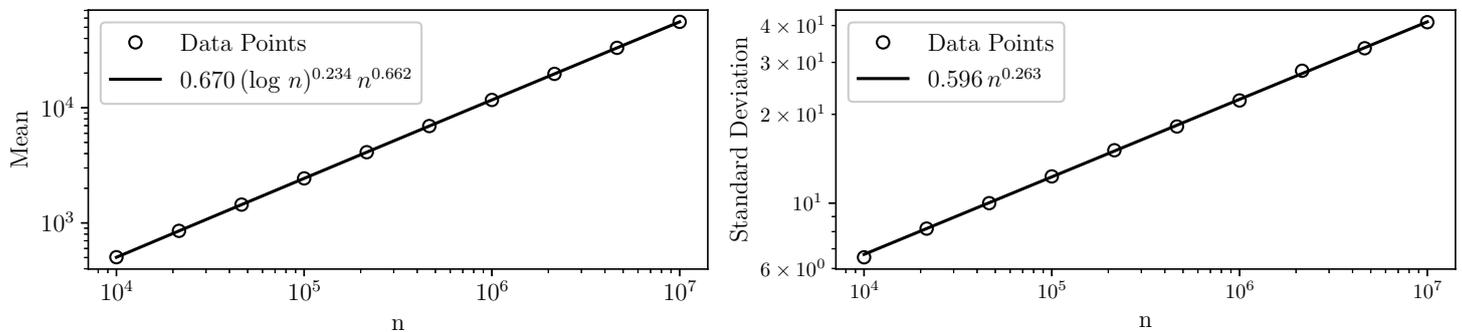

Figure 4: Log-Log plots of the best fits for the mean and standard deviation of the number of particles in model-II. The data indicate that $E(p_n^{II}) \approx a (\log n)^b n^c$ and $SD(p_n^{II}) \approx d n^f$ with 95% confidence intervals for $a \in [0.650, 0.689]$, $b \in [0.217, 0.251]$, $c \in [0.661, 0.664]$, and $d \in [0.550, 0.641]$, $f \in [0.258, 0.268]$ with corresponding MSEs 1.040 and 0.050, respectively.

|  | $n_1$ | $n_2$ | $n_3$ | $n_4$ | $n_5$ | $n_6$ | $n_7$ | $n_8$ | $n_9$ | $n_{10}$ |
|---|---|---|---|---|---|---|---|---|---|---|
| mean$(p_n^{II}) (\log n)^{-0.234} n^{-0.662}$ | 0.671 | 0.671 | 0.670 | 0.670 | 0.670 | 0.670 | 0.670 | 0.670 | 0.670 | 0.670 |
| (st. dev. $(p_n^{II})) n^{-0.263}$ | 0.583 | 0.596 | 0.594 | 0.598 | 0.600 | 0.590 | 0.591 | 0.609 | 0.593 | 0.594 |
| skewness | 0.101 | 0.087 | 0.061 | 0.042 | 0.040 | 0.099 | 0.069 | 0.092 | 0.065 | 0.003 |
| kurtosis | 3.031 | 3.051 | 2.992 | 3.038 | 2.954 | 3.099 | 3.005 | 3.015 | 3.003 | 2.943 |

Table 4: Monte Carlo simulation results for the model-II.

**Model-III** For $(u, v) \in \mathcal{K}$, let $\mathcal{D}(u, v)$ denote rectangle whose sides are parallel to $x$ and $y$ axes with vertices $(u, v), (1, v)$, $(1, 1), (u, 1)$. Initially, there are zero particles in the system. At each step, a particle appears at a uniform random point $(u, v)$ in the unit square $\mathcal{K}$; simultaneously the closest particle in $\mathcal{D}(u, v)$ to the point $(u, v)$ (if any) disappears. Let $p_n^{III}$ denote the number of particles in the system after $n$ steps.

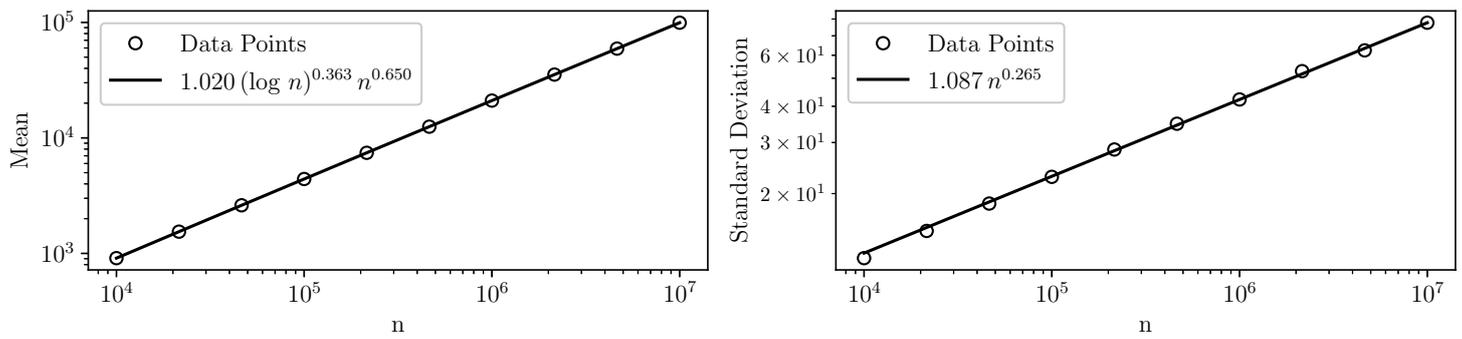

Figure 5: Log-Log plots of the best fits for the mean and standard deviation of the number of particles in model-III. The data indicate that $E(p_n^{III}) \approx a (\log n)^b n^c$ and $SD(p_n^{III}) \approx d n^f$ with 95% confidence intervals for $a \in [1.014, 1.027]$, $b \in [0.360, 0.367]$, $c \in [0.6497, 0.6502]$, and $d \in [0.978, 1.197]$, $f \in [0.258, 0.272]$ with corresponding MSEs 0.177 and 0.306, respectively.

| | $n_1$ | $n_2$ | $n_3$ | $n_4$ | $n_5$ | $n_6$ | $n_7$ | $n_8$ | $n_9$ | $n_{10}$ |
|---|---|---|---|---|---|---|---|---|---|---|
| mean$(p_n^{III}) (\log n)^{-0.363} n^{-0.650}$ | 1.020 | 1.020 | 1.020 | 1.020 | 1.020 | 1.020 | 1.020 | 1.020 | 1.020 | 1.020 |
| (st. dev. $(p_n^{III})) n^{-0.265}$ | 1.045 | 1.058 | 1.073 | 1.083 | 1.098 | 1.099 | 1.088 | 1.111 | 1.070 | 1.088 |
| skewness | 0.095 | 0.087 | 0.080 | 0.058 | 0.061 | 0.021 | 0.073 | 0.036 | 0.065 | 0.019 |
| kurtosis | 3.018 | 3.036 | 3.056 | 2.999 | 3.017 | 3.008 | 3.025 | 2.944 | 3.030 | 3.090 |

Table 5: Monte Carlo simulation results for the model-III.

## 3. Conclusion

The longest increasing subsequence for uniformly random permutations is an example of a model from the Kardar-Parisi-Zhang universality class [5]. Its study has provided a rich research program for mathematicians and physicists for the last fifty years and produced profound results in mathematics and physics, see [5,10] and references therein. Three Tracy-Widom distributions (GUE, GOE, GSE) from random matrix ensembles also appear as the limiting distributions for various subsequence problems for permutations [4,8]. For a given pattern $\tau \in S_k$, the length of the longest $\tau$-avoiding subsequence problem vastly generalizes the longest increasing subsequence problem and leads to many interesting research directions, see [1,11] for a review of the recent results. This paper will motivate further research in this direction and provide insights on the limiting distribution of the proposed models. It would be interesting to find particle processes corresponding to $\ell_n^\tau$ for general patterns $\tau$ that provide an efficient algorithm for numerical studies to understand better the limiting behaviour of the random variables.

We also studied three interacting particle systems in the plane analogous to the one-dimensional Hammersley process [2,13]. The results indicate that, in each model, the mean value and variance of the number of particles in the system after $n$ steps grows proportional to $a (\log n)^b n^c$ and $d n^f$ respectively. We also provided estimates for the corresponding exponents and coefficients. Determining the limiting distributions requires more numerical and theoretical works.

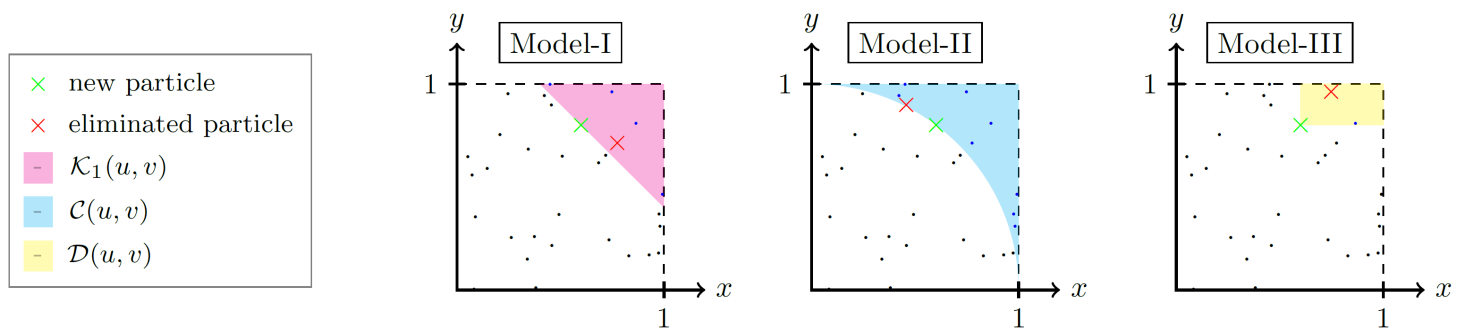

Figure 6: Illustrations of the models.

## Acknowledgements

The simulations were performed at the servers of the Department of Computer Engineering at Bilkent University. Gökhan Yıldırım was partially supported by Tubitak-Bideb through Grant No. 118C029.